# CENTRAL LIMIT THEOREM FOR STATIONARY LINEAR PROCESSES

By Magda Peligrad[1] and Sergey Utev

*University of Cincinnati and University of Nottingham*

We establish the central limit theorem for linear processes with dependent innovations including martingales and mixingale type of assumptions as defined in McLeish [*Ann. Probab.* **5** (1977) 616–621] and motivated by Gordin [*Soviet Math. Dokl.* **10** (1969) 1174–1176]. In doing so we shall preserve the generality of the coefficients, including the long range dependence case, and we shall express the variance of partial sums in a form easy to apply. Ergodicity is not required.

**1. Introduction.** Let $(\xi_i)_{i\in\mathbb{Z}}$ be a stationary sequence of random variables with $E[\xi_0^2] < \infty$ and $E[\xi_0] = 0$. Let $(a_i)_{i\in\mathbb{Z}}$ be a sequence of real numbers such that $\sum_{i\in\mathbb{Z}} a_i^2 = A < \infty$ and denote by

$$X_k = \sum_{j=-\infty}^{\infty} a_{k+j}\xi_j, \qquad S_n = \sum_{k=1}^{n} X_k,$$

(1)

$$b_{n,j} = a_{j+1} + \cdots + a_{j+n} \quad \text{and} \quad b_n^2 = \sum_{j=-\infty}^{\infty} b_{n,j}^2.$$

The so-called noncausal linear process $(X_k)_{k\in\mathbb{Z}}$ is widely used in a variety of applied fields. It is properly defined for any square summable sequence $(a_i)_{i\in\mathbb{Z}}$ if and only if the stationary sequence of innovations $(\xi_i)_{i\in\mathbb{Z}}$ has a bounded spectral density. In general, the covariances of $(X_k)_{k\in\mathbb{Z}}$ might not be summable so that the linear process might exhibit long range dependence. An important question is to describe the asymptotic properties of the variance and the asymptotic behavior of $S_n$ properly normalized. In this

Received May 2005; revised September 2005.
[1]Supported in part by a Charles Phelps Taft Memorial Fund Grant and NSA Grant H98230-05-01-0066.

paper we shall address both these questions. A simple result with very useful consequences is contained in Lemma A.3(iii). It turns out that, when the innovations have a continuous spectral density $f(x)$, the variance of $S_n$ is asymptotically proportional to $f(0)b_n^2$, up to a numerical constant. This fact suggests to further study the asymptotic distribution of $S_n/b_n$. As we shall see in this paper, if the sequence $(\xi_i)_{i\in\mathbb{Z}}$ is a martingale difference sequence or its partial sums can be approximated in a certain way by martingales, then, despite the long range dependence, $S_n/b_n$ satisfies a certain central limit theorem.

To allow for flexibility in applications, we define a stationary filtration as in [17]. We assume that $\xi_i = g(Y_j, j \leq i)$, where $(Y_i)_{i\in\mathbb{Z}}$ is an underlying stationary sequence. Denote by $\mathcal{I}$ its invariant sigma field and by $(\mathcal{F}_i)_{i\in\mathbb{Z}}$ an increasing filtration of sigma fields $\mathcal{F}_i = \sigma(Y_j, j \leq i)$. The pair $[(\mathcal{F}_i)_{i\in\mathbb{Z}}; \mathcal{I}]$ will be called a stationary filtration. For the case when for every $i$, $\xi_i = Y_j$, and $g(Y_j, j \leq i) = Y_i$, then $\mathcal{F}_i$ is simply the sigma algebra generated by $\xi_j$, $j \leq i$.

In the sequel $\|\cdot\|_2$ denotes the norm in $\mathbb{L}_2, \|X\|_2 = (E[X]^2)^{1/2}$.

We shall establish the following result:

THEOREM 1. *Let $(\xi_i)_{i\in\mathbb{Z}}$ be a stationary sequence with $E[\xi_1^2] < \infty$, $E[\xi_0] = 0$ and stationary filtration $[(\mathcal{F}_i)_{i\in\mathbb{Z}}; \mathcal{I}]$. Define $(X_k)_{k\geq 1}$, $S_n$ and $b_n$ as above and assume $b_n \to \infty$ as $n \to \infty$. Assume that*

$$\Gamma_j = \sum_{k=0}^{\infty} |E[\xi_k E(\xi_0 | \mathcal{F}_{-j})]| < \infty \quad \text{for all } j \quad \text{and}$$

(2)
$$\frac{1}{p} \sum_{j=1}^{p} \Gamma_j \to 0 \quad \text{as } p \to \infty.$$

*Then, $(\xi_i)_{i\in\mathbb{Z}}$ has a continuous spectral density $f(x)$ and there is a nonnegative random variable $\eta$ measurable with respect to $\mathcal{I}$ such that $n^{-1}E((\sum_{k=0}^{n}\xi_k)^2 | \mathcal{F}_0) \to \eta$ in $L_1$ as $n \to \infty$ and $E(\eta) = 2\pi f(0)$. In addition,*

$$\lim_{n\to\infty} \frac{\text{Var}(S_n)}{b_n^2} = 2\pi f(0) \quad \text{and}$$

(3)
$$\frac{S_n}{b_n} \Longrightarrow \sqrt{\eta} N \quad \text{in distribution as } n \to \infty,$$

*where $N$ is a standard normal variable independent of $\eta$. Moreover, if the sequence $(\xi_i)_{i\in\mathbb{Z}}$ is ergodic and condition (2) is satisfied, then the central limit theorem in (3) holds with $\eta = 2\pi f(0)$.*

The following corollary extends the projective CLT theorem of Volny [22] (which, in turn, was inspired by Heyde [11], Theorem 2) and Corollary 2 (mixingale type CLT) of Maxwell and Woodroofe [17] to dependent sequences generated by linear processes and, in addition, proves the continuity of the corresponding spectral density. This corollary also develops a result by Wu and Min [24] who considered the case of absolute summable weights.

COROLLARY 2. *Let $(\xi_i)_{i\in\mathbb{Z}}$ be a stationary sequence with $E(\xi_1^2) < \infty$, $E[\xi_0] = 0$ and stationary filtration $[(\mathcal{F}_i)_{i\in\mathbb{Z}}; \mathcal{I}]$. Consider the projection operator $P_i(Y) = E[Y|\mathcal{F}_i] - E[Y|\mathcal{F}_{i-1}]$ and assume that*

$$(4) \qquad E(\xi_0|\mathcal{F}_{-\infty}) = 0 \quad \textit{almost surely} \quad \textit{and} \quad \sum_{i=1}^{\infty} \|P_{-i}(\xi_0)\|_2 < \infty.$$

*Then, the conclusion of Theorem 1 holds. In particular, (4) is satisfied if*

$$(5) \qquad \sum_{n=1}^{\infty} n^{-1/2} \|E(\xi_n|\mathcal{F}_0)\|_2 < \infty.$$

To comment on the conditions used in our results, first we mention that assumption (2) implies that the initial sequence $(\xi_i)_{i\in\mathbb{Z}}$ satisfies the Gordin martingale approximation condition (8) defined later. Various conditions are known to be sufficient for (8), such as the original Gordin condition, $\sup_n \|E(\xi_1 + \cdots + \xi_n|\mathcal{F}_0)\|_2 < \infty$ and its modifications introduced in [11], Theorem 1, or in [9], Theorem 5.2, in [5, 7, 17, 19]. By considering telescoping sums $\xi_n = Q_n - Q_{n-1}$ with the stationary sequence $(Q_i)_{i\in\mathbb{Z}}$ having an unbounded spectral density, one can easily show that those conditions are not enough for (3). On the other hand, examples similar to those in [22], Theorem 7, show that the Gordin type conditions mentioned above, imposed to partial sums, are not necessary for (3) and (4). As a matter of fact, we shall construct an example to show that the conditions of Corollary 2 are optimal.

PROPOSITION 3. *Let $\psi_i$ be a sequence of nonnegative numbers such that $\psi_n \to 0$ as $n \to \infty$. Then, there exists a strictly stationary ergodic sequence $(\xi_i)_{i\in\mathbb{Z}}$ with unbounded spectral density such that*

$$\sum_{n=1}^{\infty} \frac{\psi_n}{n^{1/2}} \|E(\xi_n|\mathcal{F}_{-\infty}^0)\|_2 < \infty \quad \textit{and} \quad \sum_{n=1}^{\infty} \psi_n \|P_{-n}(\xi_0)\|_2 < \infty.$$

It seems that even for martingales our result is new and extends the CLT of Ibragimov [13] for linear processes with i.i.d. innovations and also the CLT of Billingsley [1] and Ibragimov [14] for stationary ergodic martingale differences to linear processes of stationary martingale differences. It also incorporates corresponding results by Heyde [11] and Hannan [10].

PROPOSITION 4. *Let $(\xi_i)_{i\in\mathbb{Z}}$ be a stationary sequence of martingale differences with finite second moment $\sigma^2$. Then (3) holds. Moreover, one can choose $\eta = E(\xi_0^2|\mathcal{I})$. In particular, if the martingale difference is ergodic, $\eta = E(\xi_0^2|\mathcal{I}) = \sigma^2$.*

The paper is organized as follows. Proofs are given in Section 2. Various examples are collected in Section 3. Among them is an application to strongly mixing structures that provides a sharp result under minimal assumptions. This section also contains the proof of Proposition 3. Finally, the Appendix gathers some technical facts about some sequences of numbers and spectral densities of stationary processes summarized as a few lemmas.

## 2. Proofs.

PROOF OF PROPOSITION 4. Denoting by $b_{n,j} = a_{j+1} + \cdots + a_{j+n}$, we express the sum $S_n/b_n = (1/b_n)\sum_{j=-\infty}^{\infty} b_{n,j}\xi_j$ and apply the central limit theorem for the triangular array of martingale differences $(b_{n,j}\xi_j/b_n)_{j\in\mathbb{Z}}$, as it was done in [20], pages 448–449, where the Lindeberg condition was established. We have only to verify the convergence condition

(6) $$\frac{1}{b_n^2}\sum_{j=-\infty}^{\infty} b_{n,j}^2 \xi_j^2 \to \eta \quad \text{in probability as } n\to\infty.$$

We start the proof by fixing a positive integer $p$ and by making small blocks of normalized sums of consecutive random variables. Define $t_{n,k} = p^{-1}\sum_{i=p(k-1)+1}^{pk} b_{n,i}^2$ and decompose the sum in (6) in the following way:

$$\frac{1}{b_n^2}\sum_{j=-\infty}^{\infty} b_{n,j}^2 \xi_j^2 = \frac{1}{b_n^2}\sum_{k=-\infty}^{\infty} pt_{n,k}\left(\frac{1}{p}\sum_{j=p(k-1)+1}^{pk} \xi_j^2\right)$$

$$+ \frac{1}{b_n^2}\sum_{k=-\infty}^{\infty}\sum_{j=p(k-1)+1}^{pk} [b_{n,j}^2 - t_{n,k}]\xi_j^2$$

$$= J_{n,1} + J_{n,2}.$$

Notice first that $\sum_k pt_{nk} = b_n^2$ and, as a consequence, by stationarity and the $\mathbb{L}_1$ ergodic theorem, the following convergence holds uniformly in $n$:

$$E|J_{n,1} - E(\xi_0^2|\mathcal{I})| \le E\left|\left(\frac{1}{p}\sum_{j=p(k-1)+1}^{pk} \xi_j^2\right) - E(\xi_0^2|\mathcal{I})\right| \to 0 \quad \text{as } p\to\infty.$$

It remains to notice that by relation (A.1) in Lemma A.1 from the Appendix it follows that

$$E|J_{n,2}| \le E[\xi_0^2]\frac{1}{b_n^2}\sum_{k=-\infty}^{\infty}\sum_{j=p(k-1)+1}^{pk} |b_{n,j}^2 - t_{n,k}| \to 0 \quad \text{as } n\to\infty.$$

☐

PROOF OF THEOREM 1. In order to prove this theorem, we shall use a blocking technique and then we shall approximate the sums of variables in blocks by martingale differences. As before, let $p$ be a fixed positive integer and denote by $I_k = \{(k-1)p+1, \ldots, kp\}$. So $I_k$'s are blocks of consecutive integers of size $p$ and $\mathbb{Z} = \bigcup_{k=-\infty}^{\infty} I_k$. We start with the following decomposition:

$$W_n := \frac{1}{b_n} \sum_{j=-\infty}^{\infty} b_{n,j} \xi_j = \frac{1}{b_n} \sum_{k=-\infty}^{\infty} \sum_{j \in I_k} b_{n,j} \xi_j.$$

With the notation $c_{n,k} = \frac{1}{p} \sum_{i \in I_k} b_{n,i}$, we further decompose $W_n$ into two terms:

$$W_n = \frac{1}{b_n} \sum_{k=-\infty}^{\infty} \sqrt{p} c_{n,k} \left( \frac{1}{\sqrt{p}} \sum_{j \in I_k} \xi_j \right) + \frac{1}{b_n} \sum_{k=-\infty}^{\infty} \sum_{j \in I_k} [b_{n,j} - c_{n,k}] \xi_j$$
$$= B_{n,1} + B_{n,2}.$$

We shall show first that $B_{n,2}$ is negligible for the convergence in distribution. Notice that by condition 2 and Lemma A.3(ii), $(\xi_i)_{i \in \mathbb{Z}}$ has a continuous spectral density and by the second inequality in part (i) of Lemma A.3, the variance of $B_{n,2}$ is bounded by

$$E(B_{n,2})^2 \leq \left( E[\xi_0^2] + 2 \sum_{k=1}^{\infty} |E(\xi_0 \xi_k)| \right) \frac{1}{b_n^2} \sum_{k=-\infty}^{\infty} \sum_{j \in I_k} [b_{n,j} - c_{n,k}]^2,$$

whence, by Lemma A.1 and taking into account condition (2), it follows that

$$E(B_{n,2})^2 \to 0 \quad \text{as } n \to \infty.$$

To analyze $B_{n,1}$, we denote the weighted sum in a block of size $p$ by

$$Y_k^{(p)} = \frac{1}{\sqrt{p}} \sum_{j \in I_k} \xi_j, \quad k \in \mathbb{Z} \quad \text{and} \quad \mathcal{G}_k = \mathcal{F}_{kp}.$$

Then, $Y_k^{(p)}$ is $\mathcal{G}_k$-measurable and define

$$Z_k^{(p)} = E(Y_k^{(p)} | \mathcal{G}_{k-1}) \quad \text{and} \quad V_k^{(p)} = Y_k^{(p)} - Z_k^{(p)}.$$

Obviously $V_k^{(p)}$ is a stationary sequence of martingale differences and $Y_k^{(p)} = Z_k^{(p)} + V_k^{(p)}$. It follows that $B_{n,1}$ can be decomposed into a linear process with stationary martingale differences innovations and another one involving $Z_k^{(p)}$.

We shall show first that the term involving $Z_k^{(p)}$ is negligible for the convergence in distribution in the sense that

$$\lim_{p\to\infty}\lim_{n\to\infty}\frac{1}{b_n^2}\left\|\sum_{k=-\infty}^{\infty}\sqrt{p}c_{n,k}(Z_k^{(p)})\right\|_2^2 = 0. \tag{7}$$

By Lemma A.1, we notice that $(1/b_n^2)\sum_{k=-\infty}^{\infty}pc_{n,k}^2 \to 1$ as $n\to\infty$ and also that the coefficients $d_{n,k} = \sqrt{p}c_{n,k}$ satisfy (A.3). Therefore, according to Lemma A.3(iii), we deduce that

$$\lim_{n\to\infty}\frac{1}{b_n^2}\left\|\sum_{k=-\infty}^{\infty}\sqrt{p}c_{n,k}(Z_k^{(p)})\right\|_2^2 = 2\pi f^{(p)}(0),$$

where $f^{(p)}(x)$ denotes the spectral density of $Z_k^{(p)}$. On the other hand, since

$$2\pi f^{(p)}(0) = \lim_{n\to\infty}\frac{1}{n}\left\|\sum_{k=1}^{n}(Z_k^{(p)})\right\|_2^2,$$

in order to establish (7), it is enough to show that

$$\lim_{p\to\infty}\sum_{k=1}^{\infty}|E(Z_1^{(p)}Z_k^{(p)})| = 0.$$

First, we observe that

$$|E(Z_1^{(p)}Z_k^{(p)})| = \frac{1}{p}\left|E\left[E\left(\sum_{i=1}^{p}\xi_i|\mathcal{F}_0\right)\sum_{j=[k-1]p+1}^{kp}\xi_j\right]\right|.$$

By the triangle inequality and condition (2), obviously

$$\sum_{k=1}^{\infty}|E(Z_1^{(p)}Z_k^{(p)})| \le 2\frac{1}{p}\sum_{i=1}^{p}\sum_{n=i}^{\infty}|E[E(\xi_i|\mathcal{F}_0)\xi_n]|$$

$$\le 2\frac{1}{p}\sum_{i=1}^{p}\Gamma_i \to 0 \quad \text{as } p\to\infty.$$

To complete the proof, we have to show that the remaining linear process involving the martingale differences satisfies the desired CLT. We shall denote by

$$X_k^{(p)} = \sum_{j=-\infty}^{\infty}\sqrt{p}c_{n,k}V_j^{(p)} \quad \text{and} \quad S_n^{(p)} = \sum_{k=1}^{n}X_k^{(p)}.$$

Notice that by Lemma A.1 and Proposition 4 it follows that, for any $p$ fixed,

$$S_n^{(p)}/b_n \to \sqrt{E([V_0^{(p)}]^2|\mathcal{I})}N \quad \text{as } n\to\infty,$$

where $N$ is a standard normal variable independent of $\mathcal{I}$. In order to complete the proof, by theorem Theorem 3.2 in [2], we have only to establish that

$$E([V_0^{(p)}]^2|\mathcal{I}) \to \eta \qquad \text{as } p \to \infty.$$

With this aim, let $T_n = \xi_1 + \cdots + \xi_n$. By applying the above decomposition and arguments to partial sums (the case $a_0 = 1$ and $a_j = 0$ for $j \geq 1$), we deduce that we have the following martingale approximation:

$$\text{(8)} \qquad \lim_{m \to \infty} \lim_{n \to \infty} \left\| n^{-1/2} T_n - \sum_{j=1}^{[n/m]} (V_j^{(m)}) \right\|_2 = 0,$$

where $[x]$ denotes the integer part of $x$, implying that Gordin's condition [8] is satisfied. Thus, by Proposition 1 in [4], there exists a nonnegative variable $\eta$ measurable with respect to $\mathcal{I}$, such that

$$E|E(p^{-1}T_p^2|\mathcal{F}_0) - \eta| \to 0 \qquad \text{as } p \to \infty.$$

It follows that $E(p^{-1}T_p^2)|\mathcal{I}) \to \eta$ and also $E(p^{-1}T_p^2) \to 2\pi f(0) = E\eta$, completing the proof of the theorem. $\square$

PROOF OF COROLLARY 2. By using a standard representation technique as in [9], by the first part of condition (4), we can write

$$\xi_k = \sum_{i=-\infty}^{k} P_i(\xi_k) \quad \text{and} \quad E(\xi_0|\mathcal{F}_{-j}) = \sum_{i=-\infty}^{-j} P_i(\xi_0).$$

By stationarity, $\|P_{-n}(\xi_0)\|_2 = \|P_{-n+k}(\xi_k)\|_2$ for any $k$. Next, $P_i(\xi_0)$ and $P_j(\xi_k)$ are uncorrelated for $i \neq j$, implying that

$$E[\xi_k E(\xi_0|\mathcal{F}_{-j})] = \sum_{i=-\infty}^{-j} E[P_i(\xi_k)P_i(\xi_0)].$$

As a consequence,

$$|E[\xi_k E(\xi_0|\mathcal{F}_{-j})]| \leq \sum_{i=-\infty}^{-j} \|P_i(\xi_k)\|_2 \|P_i(\xi_0)\|_2$$

$$= \sum_{i=-\infty}^{-j} \|P_{i-k}(\xi_0)\|_2 \|P_i(\xi_0)\|_2.$$

Therefore,

$$t_j = \sum_{k=0}^{\infty} |E[\xi_k E(\xi_0|\mathcal{F}_{-j})]| \leq \sum_{i=-\infty}^{-j} \|P_i(\xi_0)\|_2 \sum_{k=1}^{\infty} \|P_{-k}(\xi_0)\|_2,$$

whence, by (4), we derive that $\lim_{j\to\infty} t_j = 0$, that proves the validity of condition (2).

Now, we assume that (5) holds. Obviously, by the martingale convergence theorem and stationarity, $\|E(\xi_n|\mathcal{F}_0)\|_2$ is decreasing to $\|E(\xi_0|\mathcal{F}_{-\infty})\|_2$ as $n \to \infty$ and by (5), we deduce that $\|E(\xi_0|\mathcal{F}_{-\infty})\|_2 = 0$, so that the first part, of condition (4) follows. To verify its second part, we denote by $a_i := \|P_{-i}(\xi_0)\|_2 = \|P_{-i+k}(\xi_k)\|_2$ for all $k \in \mathbb{Z}$, and notice that

$$\|E(\xi_n|\mathcal{F}_0)\|_2^2 = \sum_{i=-\infty}^{0} \|P_i(\xi_n)\|_2^2 = \sum_{i=n}^{\infty} a_i^2.$$

Therefore, condition (5) and Lemma A.2 from the Appendix imply

$$\sum_{i=1}^{\infty} \|P_{-i}(\xi_0)\|_2 \leq 3 \sum_{n=1}^{\infty} n^{-1/2} \left(\sum_{i=n}^{\infty} a_i^2\right)^{1/2}$$

$$= 3 \sum_{n=1}^{\infty} n^{-1/2} \|E(\xi_n|\mathcal{F}_0)\|_2 < \infty$$

and the proof is now complete. □

### 3. Examples.

*Functionals of i.i.d. sequences.* We shall start this section by applying Corollary 2 to functionals of i.i.d. sequences. We shall see later that condition (9) required by this corollary is sharp.

COROLLARY 5. *For an i.i.d. sequence of random* $(Y_i)_{i\in\mathbb{Z}}$, *denote by* $\mathcal{F}_a^b$ *the* $\sigma$-*field generated by* $Y_k$ *with* $a \leq k \leq b$ *and define*

$$\xi_k = f(\ldots, Y_{k-1}, Y_k), \qquad k \in \mathbb{Z}.$$

*Assume that* $E(\xi_1^2) < \infty$, $E(\xi_1) = 0$ *and*

(9) $$\sum_{n=1}^{\infty} \frac{1}{\sqrt{n}} \|\xi_0 - E(\xi_0|\mathcal{F}_{-n}^0)\|_2 < \infty.$$

*Then,* (5) *is satisfied and the conclusion of Theorem 1 holds.*

PROOF. Observe that $E(\xi_0|\mathcal{F}_{-\infty}^{-n}) = E((\xi_0 - E(\xi_0|\mathcal{F}_{1-n}^0))|\mathcal{F}_{-\infty}^{-n}) + E(E(\xi_0|\mathcal{F}_{1-n}^0)|\mathcal{F}_{-\infty}^{-n})$. Now, the sigma-fields $\mathcal{F}_{1-n}^0$ and $\mathcal{F}_{-\infty}^{-n}$ are independent and so, the second term is equal almost surely to $E[E(\xi_0|\mathcal{F}_{1-n}^0)] = 0$. Therefore,

(10) $$\|E(\xi_0|\mathcal{F}_{-\infty}^{-n})\|_2 = \|E[\xi_0(\xi_0 - E(\xi_0|\mathcal{F}_{1-n}^0))]\|_2$$
$$\leq \|\xi_0\|_2 \|\xi_0 - E(\xi_0|\mathcal{F}_{1-n}^0)\|_2,$$

implying that

$$\sum_{n=1}^{\infty} \frac{1}{\sqrt{n}} \|E(\xi_n|\mathcal{F}_{-\infty}^0)\|_2 = \sum_{n=1}^{\infty} \frac{1}{\sqrt{n}} \|E(\xi_0|\mathcal{F}_{-\infty}^{-n})\|_2$$

$$\leq \|\xi_0\|_2 \sum_{n=1}^{\infty} \frac{1}{\sqrt{n}} \|\xi_0 - E(\xi_0|\mathcal{F}_{1-n}^0)\|_2 < \infty. \qquad \square$$

The following result extends Proposition 3 in [17] in the context of Bernoulli shifts (also called Raikov or Riesz–Raikov sums) and follows as an application of Corollary 2.

Let $(\varepsilon_k)_{k\in\mathbb{Z}}$ be an i.i.d. sequence with $\mathbb{P}(\varepsilon_1 = 0) = \mathbb{P}(\varepsilon_1 = 1) = 1/2$ and let

$$Y_n = \sum_{k=0}^{\infty} 2^{-k-1} \varepsilon_{n-k} \quad \text{and} \quad \xi_n = g(Y_n) - \int_0^1 g(x)\,dx,$$

where $g \in \mathbb{L}_2(0,1)$, $(0,1)$ being equipped with the Lebesgue measure.

COROLLARY 6. *For the Bernoulli shift process, if $g \in \mathbb{L}_2(0,1)$ and*

$$(11) \qquad \int_0^1 \int_0^1 [g(x)-g(y)]^2 \frac{1}{|x-y|} \left(\log\left[\log\frac{1}{|x-y|}\right]\right)^t dx\,dy < \infty$$

*for some $t > 1$, then (5) is satisfied and the conclusion of Theorem 1 holds with $\eta = 2\pi f(0)$.*

As a concrete example of a map we can take $g(x) = x^{-p}[1+\log(2/x)]^{-a} \times \sin(1/x)$, $0 < x < 1$, where either $0 \leq p < 1/2$ or $p = 1/2$ and $a > 4$. The convergence of the integral (11) is established in the same way as it was indicated in [17].

We notice that the above Corollary 5, when specified to the Bernoulli shifts, improves Theorem 19.3.1 in [15], originally established in [12, 13] and motivated by Kac [16].

PROOF OF PROPOSITION 3. We shall construct now an example to show that the conditions of Corollaries 2 and 5 are optimal. Let $(Y_i)_{i\in Z}$ be a sequence of i.i.d. random variables and assume that $Y_1$ has a standard normal distribution. As before, denote by $\mathcal{F}_a^b$ the sigma-field generated by variables $Y_k$ with $a \leq k \leq b$. Define the innovations $(\xi_i)_{i\in\mathbb{Z}}$ as a linear process

$$\xi_k = \sum_{j=-\infty}^{k} u_{k-j} Y_j,$$

where $\{u_i; i \geq 0\}$ is a sequence of nonnegative numbers to be specified. For $i < 0$, let $u_i = 0$. First, we notice that $P_{-k}(\xi_0) = u_k Y_{-k}$ and $\|P_{-k}(\xi_0)\|_2 = u_k$. Therefore,

$$\text{(12)} \quad \sum_{k=0}^{\infty} \|P_{-k}(\xi_0)\|_2 < \infty \quad \text{if and only if} \quad \sum_{k=0}^{\infty} u_k < \infty.$$

Notice that $E[\xi_0 \xi_k] = \sum_{j=0}^{\infty} u_{k+j} u_j$ and for any positive integer $j_0$, we have

$$\sum_{k=0}^{\infty} E[\xi_0 \xi_k] = \sum_{k=0}^{\infty} \sum_{j=0}^{\infty} u_{k+j} u_j > u_{j_0} \sum_{k=0}^{\infty} u_{k+j_0}.$$

So by Lemma A.3(ii), the spectral density is bounded if and only if $\sum_{j=0}^{\infty} u_j < \infty$. In particular, combining this remark with relation (12) along with the conclusion of Corollary 2, it follows that Theorem 1 holds in this example if and only if condition (2) is satisfied.

To construct $u_i$'s, without loss of generality, assume that $\psi_n \downarrow 0$. Let $n_1 = 1$ and $n_k \uparrow \infty$ be such that for $k \geq 1$, $n_{k+1} - n_k > n_{k+1}/2$ and $\psi_j \leq 1/k^2$ when $j \geq n_k$. Now, for nonnegative integers $j$, let $u_j = 1/n_{k+1}$ when $n_k \leq j < n_{k+1}$. By construction,

$$\sum_{i=1}^{\infty} u_i = \sum_{k=1}^{\infty} \sum_{n_k}^{n_{k+1}-1} 1/n_{k+1} \geq \tfrac{1}{2} \sum_{k=1}^{\infty} 1 = \infty$$

and, therefore, by the above considerations, the stationary sequence $(\xi_k)_{k \in \mathbb{Z}}$ has unbounded spectral density. By (10), it remains to show that

$$I := \sum_{j=1}^{\infty} \frac{\psi_j}{\sqrt{j}} \|(\xi_0 - E(\xi_0 | \mathcal{F}_{-1-j}^0))\|_2 < \infty.$$

Notice that since $\xi_0 - E(\xi_0 | \mathcal{F}_{-n}^0) = \sum_{i=-\infty}^{-n-1} u_{-i} Y_i$, and also, since, for $j \geq n_k$,

$$\sum_{i=j}^{\infty} u_i^2 \leq \sum_{i=k}^{\infty} \sum_{j=n_i}^{n_{i+1}-1} u_j^2 \leq \sum_{i=k}^{\infty} (1/n_{i+1}^2) n_{i+1} \leq c_1 \frac{1}{n_{k+1}},$$

we derive the following estimate:

$$I = \sum_{k=1}^{\infty} \sum_{j=n_k}^{n_{k+1}-1} \frac{\psi_j}{\sqrt{j}} \left[ \sum_{i=j}^{\infty} u_i^2 \right]^{1/2}$$

$$\leq \sum_{k=1}^{\infty} \frac{1}{k^2} \sum_{j=n_k}^{n_{k+1}-1} \left[ \frac{1}{j} \sum_{i=j}^{\infty} u_i^2 \right]^{1/2}$$

$$\leq \sum_{k=1}^{\infty} \frac{1}{k^2} \frac{c_1}{\sqrt{n_{k+1}}} \sum_{j=1}^{n_{k+1}} \frac{1}{\sqrt{j}} < \infty. \qquad \square$$

*Mixingales.* We are going to apply Theorem 1 to mixingales and strongly mixing sequences. For a stationary sequence of random variables $(\xi_k)_{k \in \mathbb{Z}}$, we define $\mathcal{F}_m^n$ the sigma-field generated by $\xi_i$ with indices $m \leq i \leq n$ and the sequences of coefficients $\alpha(n)$:

$$\alpha(n) = \alpha(\mathcal{F}_{-\infty}^0, \mathcal{F}_n^\infty) = \sup\{|\mathbb{P}(A \cap B) - \mathbb{P}(A)\mathbb{P}(B)|; A \in \mathcal{F}_{-\infty}^0, B \in \mathcal{F}_n^\infty\}.$$

We say that the strictly stationary sequence is strongly mixing if $\alpha(n) \to 0$ as $n \to \infty$. Various examples of mixing sequences can be found in books by Rio [21] and Bradley [3], along with counterexamples showing that the conditions we use in the next corollary are sharp for central limit theorem even for partial sum processes. In the next corollary we shall use a weaker form of the strongly mixing coefficient, a mixingale type condition, where $\mathcal{F}_n^\infty$ is replaced by the sigma-field generated by $\xi_n$, namely, $\bar{\alpha}(n) = \alpha(\mathcal{F}_{-\infty}^0, \mathcal{F}_n^n)$.

COROLLARY 7. *Assume that the innovations $(\xi_k, k \in \mathbb{Z})$ form a stationary sequence of centered random variables with finite second moment and such that*

$$\text{(13)} \qquad \sum_{k=1}^\infty \int_0^{\bar{\alpha}(k)} Q^2(u)\, du < \infty,$$

*where $Q$ denotes the cadlag inverse of the function $t \to P(|\xi_0| > t)$.*

*Then the conclusion of Theorem 1 holds. Moreover, with $\bar{\alpha}(k)$ being replaced by $\alpha(k)$, the sequence $\xi_k$ is ergodic and $\eta$ is a constant $\eta = 2\pi f(0)$, where $f(x)$ is the continuous spectral density of the innovations.*

PROOF. According to Theorem 1, it is enough to establish the validity of the condition (2). We notice that, by Rio's (1993) covariance inequality (see also [21], Chapter 4), we have

$$|E(\xi_k E(\xi_0|\mathcal{F}_{-j}))| \leq \int_0^{\bar{\alpha}(k+j)} Q^2(u)\, du,$$

that proves that condition (2) of Theorem 1 holds, since

$$\sum_{k=0}^\infty |E(\xi_k E(\xi_0|\mathcal{F}_{-j}))| \leq \sum_{i=j}^\infty \int_0^{\bar{\alpha}(i)} Q^2(u)\, du \to 0 \qquad \text{as } j \to \infty. \qquad \square$$

In comparison with Peligrad and Utev [20], Corollary 7 provides explicit normalizing constants.

To make condition (13) more transparent, we mention that it is implied by the couple of conditions (as it was derived in [6])

$$E|X_0|^t < \infty \quad \text{and} \quad \sum_{k=1}^\infty k^{2/(t-2)} \tilde{\alpha}(k) < \infty \qquad \text{where } t > 2.$$

## APPENDIX

*Facts about sequences.*

LEMMA A.1. *Let $b_{n,j} = a_{j+1} + \cdots + a_{j+n}$, for $j \in \mathbb{Z}$ and $n \in \mathbb{N}$. Assume that*

$$b_n^2 = \sum_{j=-\infty}^{\infty} b_{n,j}^2 \to \infty \quad \text{and} \quad \sum_{j \in \mathbb{Z}} a_j^2 < \infty.$$

*Then,*

(A.1) $\quad \dfrac{1}{b_n^2} \displaystyle\sum_{j=-\infty}^{\infty} |b_{n,j} - b_{n,j-1}|^2 \to 0 \quad \text{and} \quad \dfrac{1}{b_n^2} \displaystyle\sum_{j=-\infty}^{\infty} |b_{n,j}^2 - b_{n,j-1}^2| \to 0.$

*More generally, let $p$ be a positive integer. Starting with zero (in two directions), we denote blocks of consecutive integers of size $p$ by $I_k$. For each $k$, define averages of the $b_{n,i}$ in $I_k$ by $c_{n,k} = \frac{1}{p} \sum_{i \in I_k} b_{n,i}$. Then, as $n \to \infty$,*

(A.2) $\quad \dfrac{1}{b_n^2} \displaystyle\sum_{k \in \mathbb{Z}} \sum_{j \in I_k} |b_{n,j} - c_{n,k}|^2 \to 0 \quad \text{and} \quad \dfrac{1}{b_n^2} \displaystyle\sum_{k \in \mathbb{Z}} \sum_{j \in I_k} |b_{n,j}^2 - c_{n,k}^2| \to 0.$

PROOF. To simplify the writing, let us denote by $(b_n')^2 = \sum_{j=-\infty}^{\infty} |b_{n,j} - b_{n,j-1}|^2$. The validity of the first part of relation (A.1) is straightforward from the following observation:

$$(b_n')^2 \leq \sum_{j=-\infty}^{\infty} |a_j - a_{n+j+1}|^2 \leq 4 \sum_{j=-\infty}^{\infty} a_j^2,$$

implying that $\lim_{n \to \infty} (b_n'/b_n)^2 = 0$. The second part easily follows by applying Hölder inequality:

$$\sum_{j=-\infty}^{\infty} |b_{n,j}^2 - b_{n,j-1}^2| = \sum_{j=-\infty}^{\infty} |b_{n,j} - b_{n,j-1}| * |b_{n,j} + b_{n,j-1}| \leq C b_n' b_n.$$

The proof of (A.2) is similar by taking into account that $p$ is a fixed positive integer and for any pair of indexes $i, l \in I_k$, we have

$$|b_{n,i} - b_{n,l}|^2 \leq p \sum_{j \in I_k} |b_{n,j} - b_{n,j-1}|^2. \qquad \square$$

LEMMA A.2. *Suppose that $(a_j)_{j \in \mathbb{N}}$ is a sequence of nonnegative numbers and $\psi_n$ is a nonincreasing sequence of nonnegative numbers. Then,*

$$\sum_{n=1}^{\infty} a_n \psi_n \leq 3 \sum_{n=1}^{\infty} n^{-1/2} \psi_n \left( \sum_{k=n}^{\infty} a_k^2 \right)^{1/2}.$$

PROOF. The proof involves an application of the inequality in [23], contained in his Lemma 1 to $g_n = a_n \psi_n$, $n = 1, 2, \ldots$. We obtain

$$\sum_{n=1}^{\infty} a_n \psi_n \leq 3 \sum_{n=1}^{\infty} n^{-1/2} \left( \sum_{k=n}^{\infty} \psi_k^2 a_k^2 \right)^{1/2}$$

$$\leq 3 \sum_{n=1}^{\infty} n^{-1/2} \psi_n \left( \sum_{k=n}^{\infty} a_k^2 \right)^{1/2},$$

where at the last step we have used the fact that the sequence $\psi_n$ is nondecreasing. $\square$

*Facts about spectral densities.* In the following lemma we combine a few facts about spectral densities, covariances, behavior of variances of sums and their relationships. The first two points are known and can be found in books by Bradley [3].

LEMMA A.3. *Let $(\xi_i)_{i \in \mathbb{Z}}$ be a stationary sequence of real valued variables with $E[\xi_0] = 0$ and finite second moment. Let $F$ denotes the spectral measure and $f$ denotes its spectral density (if exists), that is,*

$$E[\xi_0 \xi_k] = \int_{-\pi}^{\pi} e^{-ikt} \, dF(t) = \int_{-\pi}^{\pi} e^{-ikt} f(t) \, dt.$$

(i) *For any positive integer $n$ and any real numbers $a_1, \ldots, a_n$,*

$$E\left( \sum_{k=1}^{n} a_k \xi_k \right)^2 = \int_{-\pi}^{\pi} \left| \sum_{k=1}^{n} a_k e^{ikt} \right|^2 f(t) \, dt \leq 2\pi \|f\|_{\infty} \sum_{k=1}^{n} a_k^2$$

$$\leq \left( E[\xi_0^2] + 2 \sum_{k=1}^{\infty} |E(\xi_0 \xi_k)| \right) \sum_{k=1}^{n} a_k^2.$$

(ii) *Assume* (B): $\sum_{k=1}^{\infty} |E(\xi_0 \xi_k)| < \infty$. *Then, $f$ is continuous. Moreover, if $E[\xi_k \xi_0] \geq 0$ for all $k$, then the spectral density is bounded if and only if relation* (B) *is satisfied.*

(iii) *Assume that the spectral density $f$ is continuous, and let $d^{(n)} = (d_{n,j})_{j \in \mathbb{Z}}$ be a double array of real numbers with $d_n^2 = \sum_{j \in \mathbb{Z}} d_{n,j}^2 < \infty$ that satisfies the condition*

(A.3) $$\frac{1}{d_n^2} \sum_{j=-\infty}^{\infty} |d_{n,j} - d_{n,j-1}|^2 \to 0.$$

*Then,*

(A.4) $$\lim_{n \to \infty} \frac{1}{d_n^2} E\left( \sum_{j=1}^{n} d_{n,j} \xi_j \right)^2 = 2\pi f(0).$$

PROOF. (iii) Fix $\varepsilon > 0$. By the Stone–Weierstrass theorem, there exists a trigonometric polynomial $P_m(t) = \sum_{k=-m}^{m} c_k e^{itk}$ such that $\sup_{t \in [-\pi,\pi]} |f(t) - P_m(t)| \leq \varepsilon$. In particular,

$$\left| f(0) - \sum_{k=-m}^{m} c_k \right| \leq \varepsilon. \tag{A.5}$$

Whence, by (i),

$$\frac{1}{d_n^2} E\left( \sum_{j \in \mathbb{Z}} d_{n,j} \xi_j \right)^2 = \frac{1}{d_n^2} \int_{-\pi}^{\pi} \left| \sum_{j \in \mathbb{Z}} d_{n,j} e^{itj} \right|^2 f(t)\, dt$$

$$= O(\varepsilon) + \frac{1}{d_n^2} \int_{-\pi}^{\pi} \left| \sum_{j \in \mathbb{Z}} d_{n,j} e^{itj} \right|^2 P_m(t)\, dt. \tag{A.6}$$

With the notation $A_{n,k} := d_n^{-2} \sum_{j \in \mathbb{Z}} d_{n,j} d_{n,j+k}$, we have

$$\frac{1}{d_n^2} \int_{-\pi}^{\pi} \left| \sum_{j \in \mathbb{Z}} d_{n,j} e^{itj} \right|^2 P_m(t)\, dt = 2\pi \sum_{k=-m}^{m} c_k A_{n,k}. \tag{A.7}$$

By (A.3) and similar arguments as in the proof of Lemma A.1, we can see that $A_{n,k}$ can be easily approximated by $d_n^{-2} \sum_{j \in \mathbb{Z}} d_{n,j}^2$ and, as a consequence, for any $k$ fixed, $A_{n,k}$ approaches 1 as $n \to \infty$, implying that $\lim_{n \to \infty} 2\pi \sum_{k=-m}^{m} c_k A_{n,k} = 2\pi \sum_{k=-m}^{m} c_k$. We have now only to combine this convergence with (A.5) and (A.6) to complete the proof of the statement. □

**Acknowledgment.** The authors thank the anonymous referee for helpful comments.

## REFERENCES

[1] BILLINGSLEY, P. (1961). The Lindeberg–Lévy theorem for martingales. *Proc. Amer. Math. Soc.* **12** 788–792. MR0126871
[2] BILLINGSLEY, P. (1999). *Convergence of Probability Measures*, 2nd ed. Wiley, New York. MR1700749
[3] BRADLEY, R. C. (2002, 2003). Introduction to strong mixing conditions. **1**, **2**. Technical report, Dept. Mathematics, Indiana Univ., Bloomington. Custom Publishing of I.U., Bloomington.
[4] DEDECKER, J. and MERLEVÈDE, F. (2002). Necessary and sufficient conditions for the conditional central limit theorem. *Ann. Probab.* **30** 1044–1081. MR1920101
[5] DEDECKER, J. and RIO, E. (2000). On the functional central limit theorem for stationary processes. *Ann. Inst. H. Poincaré Probab. Statist.* **36** 1–34. MR1743095
[6] DOUKHAN, P., MASSART, P. and RIO, E. (1994). The functional central limit theorem for strongly mixing processes. *Ann. Inst. H. Poincaré Probab. Statist.* **30** 63–82. MR1262892

[7] DURR, D. and GOLDSTEIN, S. (1986). Remarks on the central limit theorem for weakly dependent random variables. *Lecture Notes in Math.* **1158** 104–118. Springer, Berlin. MR0838560
[8] GORDIN, M. I. (1969). The central limit theorem for stationary processes. *Soviet Math. Dokl.* **10** 1174–1176. MR0251785
[9] HALL, P. and HEYDE, C. C. (1980). *Martingale Limit Theory and Its Application.* Academic Press, New York. MR0624435
[10] HANNAN, E. J. (1979). The central limit theorem for time series regression. *Stochastic Process. Appl.* **9** 281–289. MR0562049
[11] HEYDE, C. C. (1974). On the central limit theorem for stationary processes. *Z. Wahrsch. Verw. Gebiete* **30** 315–320. MR0372955
[12] IBRAGIMOV, I. A. (1960). On asymptotic distribution of values of certain sums. *Vestnik Leningrad. Univ.* **15** 55–69. (In Russian.) MR0120679
[13] IBRAGIMOV, I. A. (1962). Some limit theorems for stationary processes. *Theory Probab. Appl.* **7** 349–382. MR0148125
[14] IBRAGIMOV, I. A. (1963). A central limit theorem for a class of dependent random variables. *Theory Probab. Appl.* **8** 83–89. MR0151997
[15] IBRAGIMOV, I. A. and LINNIK, YU. V. (1971). *Independent and Stationary Sequences of Random Variables*. Wolters, Groningen. MR481287
[16] KAC, M. (1946). On the distribution of values of sums of the type $\sum f(2^k t)$. *Ann. of Math. (2)* **47** 33–49. MR0015548
[17] MAXWELL, M. and WOODROOFE, M. (2000). Central limit theorems for additive functionals of Markov chains. *Ann. Probab.* **28** 713–724. MR1782272
[18] MCLEISH, D. L. (1977). On the invariance principle for nonstationary mixingales. *Ann. Probab.* **5** 616–621. MR0445583
[19] PELIGRAD, M. (1981). An invariance principle for dependent random variables. *Z. Wahrsch. Verw. Gebiete* **57** 495–507. MR0631373
[20] PELIGRAD, M. and UTEV, S. (1997). Central limit theorem for linear processes. *Ann. Probab.* **25** 443–456. MR1428516
[21] RIO, E. (2000). *Théorie asymptotique des processus aléatoires faiblement dépendants. Mathém. Appl. SMAI* **31**. Springer, Berlin. MR2117923
[22] VOLNÝ, D. (1993). Approximating martingales and the central limit theorem for strictly stationary processes. *Stochastic Process. Appl.* **44** 41–74. MR1198662
[23] WU, W. B. (2002). Central limit theorems for functionals of linear processes and their applications. *Statist. Sinica* **12** 635–649. MR1902729
[24] WU, W. B. and MIN, W. (2005). On linear processes with dependent innovations. *Stochastic Process. Appl.* **115** 939–958. MR2138809

DEPARTMENT OF MATHEMATICAL SCIENCES
UNIVERSITY OF CINCINNATI
PO BOX 210025
CINCINNATI, OHIO 45221-0025
USA
E-MAIL: peligrm@math.uc.edu

SCHOOL OF MATHEMATICAL SCIENCES
UNIVERSITY OF NOTTINGHAM
NOTTINGHAM, NG7 2RD
ENGLAND
E-MAIL: pmzsu@gwmail.nottingham.ac.uk